\documentclass[12pt]{article}
\usepackage{stmaryrd,bbm,epsfig}
\usepackage{wrapfig,epic}
\usepackage{amsmath}
\usepackage{amssymb}
\usepackage{amscd}
\usepackage{rotating}
\usepackage{color}
\newcommand{\blue}[1]{\textcolor{blue}{#1}}

\DeclareMathAlphabet{\mathpzc}{OT1}{pzc}{m}{it}

\begin{document}
\title{Singular Algebraic Equations with Empirical Data}

\author{Zhonggang Zeng \thanks{Department of Mathematics,
Northeastern Illinois University, Chicago, Illinois 60625, USA.
~~email:~{\tt zzeng@neiu.edu}. 
~Research is supported in part by NSF under grant DMS-1620337.}}

\newcommand{\qed}{\hfill $\blacksquare$}
\newtheorem{example}{Example}
\newtheorem{problem}{Problem}[section]
\newtheorem{theorem}{Theorem}
\newtheorem{lemma}{Lemma}

\newcommand{\la}{\lambda}
\newcommand{\CC}{\mathbbm{C}}
\newcommand{\blb}{\big[\,}
\newcommand{\brb}{\, \big]}
\newcommand{\h}{{{\mbox{\tiny $\mathsf{H}$}}}}
\newcommand{\cB}{{\mathcal{B}}}
\newcommand{\nullity}[1]{\mathpzc{nullity}\big(\,#1\,\big)}
\newcommand{\eps}{\varepsilon}
\newcommand{\dl}{\delta}
\newcommand{\cK}{\mathpzc{Kernel}}
\newcommand{\cR}{\mathpzc{Range}}
\newcommand{\cC}{\mathcal{C}}
\newcommand{\cP}{\mathcal{P}}
\newcommand{\cF}{\mathcal{F}}
\newcommand{\cS}{\mathcal{S}}
\newcommand{\cT}{\mathcal{T}}
\newcommand{\cU}{\mathcal{U}}
\newcommand{\cV}{\mathcal{V}}
\newcommand{\cW}{\mathcal{W}}
\newcommand{\cM}{\mathcal{M}}
\newcommand{\cN}{\mathcal{N}}
\newcommand{\bdu}{\mathbf{u}}
\newcommand{\bdv}{\mathbf{v}}
\newcommand{\bdw}{\mathbf{w}}
\newcommand{\bdx}{\mathbf{x}}
\newcommand{\bdy}{\mathbf{y}}
\newcommand{\bdz}{\mathbf{z}}
\newcommand{\bda}{\mathbf{a}}
\newcommand{\bdb}{\mathbf{b}}
\newcommand{\bdc}{\mathbf{c}}
\newcommand{\bdd}{\mathbf{d}}
\newcommand{\bde}{\mathbf{e}}
\newcommand{\bdf}{\mathbf{f}}
\newcommand{\bdg}{\mathbf{g}}
\newcommand{\bdh}{\mathbf{h}}
\newcommand{\bdi}{\mathbf{i}}
\newcommand{\bdj}{\mathbf{j}}
\newcommand{\bdm}{\mathbf{m}}
\newcommand{\bdn}{\mathbf{n}}
\newcommand{\bdk}{\mathbf{k}}
\newcommand{\bdl}{\mathbf{l}}
\newcommand{\bdr}{\mathbf{r}}
\newcommand{\bdp}{\mathbf{p}}
\newcommand{\bdq}{\mathbf{q}}
\newcommand{\bds}{\mathbf{s}}
\newcommand{\bdo}{\mathbf{0}}
\newcommand{\bdF}{\mathbf{F}}
\newcommand{\dm}{\dim}
\newcommand{\dg}{\deg}
\newcommand{\codim}{\mathrm{codim}}
\newcommand{\al}{\alpha}
\newcommand{\bt}{\beta}
\newcommand{\gm}{\gamma}
\newcommand{\Dl}{\Delta}
\newcommand{\sg}{\sigma}
\newcommand{\rank}[1]{\mathpzc{rank}\big(\,#1\,\big)}
\newcommand{\ranka}[2]{\mathpzc{rank}_{#1}\left(\,#2\,\right)}
\newcommand{\GCD}{\mathrm{gcd}}
\newcommand{\dist}[1]{\mathrm{dist}\big(\,#1\,\big)}
\newcommand{\eig}{\mathpzc{eig}}
\newcommand{\rkr}{{\mbox{\scriptsize rank-$r$}}}
\newcommand{\rk}[1]{{\mbox{\scriptsize rank-{#1}}}}
\newcommand{\spn}{\mathpzc{span}}

\newcommand{\mapform}[5]{\begin{array}{ccrcl}
#1 & : & #2 & \longrightarrow & #3 \\
& & #4 & \longmapsto & #5
\end{array}}

\maketitle


\begin{abstract}
Singular equations with rank-deficient Jacobians arise frequently in 
algebraic computing applications.
As shown in case studies in this paper, direct and intuitive modeling of 
algebraic problems often results in nonisolated singular solutions.
The challenges become formidable when the problems need to be solved
from empirical data of limited accuracy.
A newly discovered low-rank Newton's iteration emerges as an effective
regularization mechanism that enables solving singular equations accurately
with an error bound in the same order as the data error.
This paper elaborates applications of new methods on solving singular
algebraic equations such as singular linear systems, polynomial GCD and 
factorizations as well as matrix defective eigenvalue problems.
\end{abstract}

\vspace{-4mm}
\section{Introduction}

\vspace{-4mm}
Algebraic equations can have singular solutions at which the Jacobians
are rank-deficient.
Those singular solutions can be isolated or in a form of varieties of
positive dimensions.
Such singular equations pose formidable challenges in scientific computing 
especially when the data are given from measurement with limited
accuracy or processed with necessary round-off.
Direct attempt of solving singular equations from empirical data may not 
achieve accurate solutions since the solutions can be altered
substantially or even disappear.
Common iterative methods such as Newton's iteration are not guaranteed to 
converge at singular solutions.
Those difficulties are well-documented in the literature such as
\cite{BSHW13,Gri85,Kel81,DecKelKel,victorpan97,ypma}.
Even for linear equations whose singular solutions are elementary
in linar algebra, the textbook advice is still to avoid solving singular 
equations with any perturbation \cite[pp 217-218]{Meyer}.
Theories and computational methodologies appear to be inadequte on 
singular algebraic equations particularly with empirical data.

Newton's iteration is tremendously effective in solving regular equations.
However, its textbook formulation is only a special case and finding 
nonsingular solutions is merely a fraction of its capabilities.
Extending Newton's iteration to solving singular equations has been studied
in many works over the years such as 
\cite{Ben66,Chu83,DedKim,NasChe93}.
By a simple modification, a low-rank Newton's iteration emerges 
as an effective method in solving singular equations for nonisolated solutions
\cite{ZengNewton} and maintains quadratic convergence. 
More importantly, it serves as a regularization mechanism so that singular 
solutions can be solved accurately from perturbed data and the solution 
accuracy is bounded by a multiple the data error.

A generic class of singular equations that possess {\em semiregular} solutions
are emphasized in this paper as opposed to {\em ultrasingular} ones, and we 
shall elaborate the low-rank Newton's iteration on such equations.
We shall establish semiregularity of some fundamental singular equations 
such as polynomial GCD/factorization and defective eigenvalue 
equations as applications of the low-rank Newton's iteration and demonstrate
its effectiveness in their accurate solutions.
We shall also briefly elaborate experimental results of solving 
ultrasingular equations using the depth-deflation method.

We restrict our elaboration to solving singular equations as zero-finding 
for holomorphic mappings in complex domains.
The same theories and computing methods apply to singular zeros of 
real mappings that are twice continuously differentiable. 

\vspace{-4mm}
\section{Preliminaries}

\vspace{-4mm}
The space of $n$-dimensional vectors of complex numbers is denoted by 
$\CC^n$ with the Euclidean norm $\|\cdot\|_2$.
Matrices of $m\times n$ form the vector space $\CC^{m\times n}$ with
the Frobenius norm $\|\cdot\|_{_F}$.
Matrices are denoted by upper-case letters with $(\cdot)^\h$ being the 
Hermitian transpose of any matrix $(\cdot)$.
A zero matrix is denoted by $O$ whose the sizes are derived from the context.

Finite-dimensional normed vector spaces are denoted by, say $\cV$, $\cW$ etc, 
in which vectors are denoted by boldface lower-case
letters with $\bdo$ being the zero vector.
For any vector $\bdv$, the norm $\|\bdv\|$ is understood as {\em the} norm 
in the space where $\bdv$ belongs.
For any linear map $L : \cV\rightarrow\cW$, its norm is the 
operator norm 
\[ \|L\| := \max_{\bdv\in\cV, \|\bdv\|=1}\|L(\bdv)\|
\]
derived from the norms of its domain $\cV$ and codomain $\cW$.
A vector space $\cV$ can be isomorphic to $\CC^n$ where $n = 
\dm(\cV)$, the dimension of $\cV$.
Throughout this paper, the norm of a product space
$\cV\times\cW$ is 
\[ \|(\bdv,\bdw)\| := \sqrt{\|\bdv\|^2+\|\bdw\|^2} 
~~~\mbox{for}~~~ (\bdv,\bdw)\in\cV\times\cW.
\]
For any linear map $L$, the notations $\cR(L)$, $\cK(L)$, $\rank{L}$ 
and $\nullity{L}$ represent the range, kernel, rank and nullity of $L$
respectively.

For a holomorphic mapping $F : \Omega\subset\CC^n\rightarrow\CC^m$,
we can designate a variable name, say $\bdz$, and denote $F$ as
$\bdz\mapsto F(\bdz)$.
Then the Jacobian of $F$ at any $\bdz_0\in\Omega$ is the matrix denoted by 
$F_\bdz(\bdz_0)$.
Let $\cV$ and $\cW$ be normed vector spaces isomorphic
to $\CC^n$ and $\CC^m$ respectively via isomorphisms 
$\psi_{_\cV}:\cV\rightarrow\CC^n$ and $\psi_{_\cW}:\cW\rightarrow\CC^m$. 
Assume $\bdv\mapsto\bdg(\bdv)$ is a mapping from an open subset 
$\Sigma$ of $\cV$ to $\cW$ with a representation $\bdz\mapsto G(\bdz)$ where
$G:\psi_{_\cV}(\Sigma)\subset\CC^n\rightarrow
\CC^m$ such that $\bdg = \psi_{_\cW}^{-1}\circ G\circ \psi_{_\cV}$
that makes the following diagram commute
\[
\begin{CD} \Sigma\subset\cV @>{\bdg}>> \cW \\
@V{\psi_{_\cV}    }VV  @AA{ \psi_{_\cW}^{-1}  }A \\
\psi_{_\cV}(\Sigma)\subset\CC^n @>G>> \CC^m.
\end{CD}
\]
We say $\bdg$ is holomorphic in $\Sigma$ if $G$ is holomorphic in 
$\psi_{_\cV}(\Sigma)$. 
The {\em Jacobian} of $\bdg$ at any $\bdv_0\in\Sigma$ is defined as the 
linear map $\bdg_\bdv(\bdv_0)$ in the form of 
\begin{equation}\label{gjac}
\mapform{\bdg_\bdv(\bdv_0)}{\cV}{\cW}{\bdv}{
\psi_{_\cW}^{-1}\circ G_\bdz(\bdz_0)\circ\psi_{_\cV}(\bdv)}
\end{equation}
where $\bdz_0 = \psi_{_\cV}(\bdv_0)$.
The Jacobian $\bdg_\bdv(\bdv_0)$ as a linear map is invariant under change 
of bases.
Let $G_\bdz(\bdz_0)^\dagger$ be the Moore-Penrose inverse of the Jacobian 
matrix $G_\bdz(\bdz_0)$.
If we further assume the isomorphisms $\psi_{_\cV}$ and $\psi_{_\cW}$ are 
isometric, namely
\[ \|\psi_{_\cV}(\bdv)\|_2 =\|\bdv\| ~~~\mbox{and}~~~
\|\psi_{_\cW}(\bdw)\|_2 = \|\bdw\| ~~~\mbox{for all}~~~
\bdv\in\cV ~~\mbox{and}~~ \bdw\in\cW, 
\]
then $\bdg_\bdv(\bdv_0)^\dagger$ is well-defined as 
\[ \bdg_\bdv(\bdv_0)^\dagger  = \psi_{_\cV}^{-1}\circ G_\bdz(\bdz_0)^\dagger
\circ\psi_{_\cW}
\]
that is invariant under isometric isomorphisms.

For any matrix $A\in\CC^{m\times n}$, we denote $A_\rkr$ as the {\em rank-$r$
projection} of $A$. 
Namely $A_\rkr$ is the rank-$r$ matrix with the smallest distance
$\|A_\rkr-A\|_{_F}$ to $A$.
The rank-$r$ projection is also called the rank-$r$ approximation and
rank-$r$ truncated SVD in the literature.
For a holomorphic mapping $\bdv\mapsto\bdg(\bdv)$ with its Jacobian
$\bdg_\bdv(\bdv_0)$ defined in (\ref{gjac}), its {\em rank-$r$ projection}
$\bdg_\bdv(\bdv_0)_\rkr$ is defined as the linear map
\[\bdg_\bdv(\bdv_0)_\rkr \,:\, \bdv\longmapsto
\psi_{_\cW}^{-1}\circ G_\bdz(\bdz_0)_\rkr\circ\psi_{_\cV}(\bdv)
\]
where $\bdz_0=\psi_{_\cV}(\bdv_0)$ and $G_\bdz(\bdz_0)_\rkr$ is the
rank-$r$ projection of the Jacobian matrix $G_\bdz(\bdz_0)$.
The notation $\bdg_\bdv(\bdv_0)_\rkr^\dagger := 
 \big(\bdg_\bdv(\bdv_0)_\rkr\big)^\dagger$.

For multivariate mappings, say $(\bdu,\bdv,\bdw)\mapsto\bdf(\bdu,\bdv,\bdw)$, 
its Jacobian at $(\bdu_0,\bdv_0,\bdw_0)$ is denoted by 
$\bdf_{\bdu \bdv \bdw}(\bdu_0,\bdv_0,\bdw_0)$, and the notation such as
$\bdf_{\bdu \bdw}(\bdu_0,\bdv_0,\bdw_0)$ denotes
the partial Jacobian with respect to $(\bdu,\bdw)$.

\vspace{-4mm}
\section{Semiregular and ultrasingular zeros}

\vspace{-4mm}
An equation $\bdf(\bdx) = \bdo$ is {\em singular} if the Jacobian
$\bdf_\bdx(\bdx_*)$ is rank-deficient so that 
$\nullity{\bdf_\bdx(\bdx_*)}>0$ at the desired solution $\bdx_*$.
A solution $\bdx_*$ is {\em isolated} if there is an open neighborhood $\Dl$
such that $\Dl\cap\bdf^{-1}(\bdo)=\{\bdx_*\}$.
The Jacobian $\bdf_\bdx(\bdx_*)$ with nullity zero always implies 
$\bdx_*$ is isolated, and $\bdx_*$ is a {\em regular} zero.
A nonisolated solution is singular and may be a point on a curve, a surface 
etc.

For a holomorphic mapping $\bdf : \Omega\subset\cV\rightarrow\cW$,
we say the {\em dimension} of its zero $\bdx_*$ is $k$
if there is an open neighborhood $\Dl\subset\Omega$ of $\bdx_*$ in 
$\cV$ such that $\Dl\cap\bdf^{-1}(\bdo) = \phi(\Lambda)$ where 
$\bdz\mapsto\phi(\bdz)$ is a holomorphic injective mapping defined 
in a connected open set $\Lambda$ in $\CC^k$ for $k > 0$ ~with 
$\phi(\bdz_*) = \bdx_*$ and $\rank{\phi_\bdz(\bdz_*)} = k$.
As a special case, an isolated zero is of dimension $0$.
A singular zero $\bdx_*$ of $\bdf$ is said to be {\em semiregular} if its 
dimension is identical to $\nullity{\bdf_\bdx(\bdx_*)}$. 
A zero is {\em ultrasingular} if it is not semiregular.
We say an equation is semiregular or ultrasingular if the intended solutions
are.

The identity $\bdf(\phi(\bdz))\equiv\bdo$ for $\bdz$ in the domain of $\phi$
implies $\bdf_\bdx(\bdx_*)\phi_\bdz(\bdz_*)$ is a zero mapping and thus
$\nullity{\bdf_\bdx(\bdx_*)}\ge k$ since $\phi_\bdz(\bdz_*)$ is of rank $k$.
A simple approach to establish semiregularity of a $k$-dimensional solution
is to append a linear mapping $L$ to $\bdf$ from the same domain as $\bdf$ 
to a codomain of dimension $k$. 
A $k$-dimensional zero $\bdx_*$ is semiregular if Jacobian of the stacked 
mapping $\bdx\mapsto \big(\bdf(\bdx),L(\bdx)\big)$ is injective so
that $\nullity{\bdf_\bdx(\bdx_*)}\le k$ and must equal to $k$.
We shall apply this technique repeatedly in the sample applications.

Because tiny perturbations can only reduce nullities, semiregularity 
$\nullity{\bdf_\bdx(\bdx_*)} = k$ is generic among singular solutions. 
Extra singularity is required to increase the nullity further and moves away
from semiregularity.
As a special case, a regular zero is semeregular with dimension 0.

Semiregular equations enjoy the {\em stationary point property}
\cite[Lemma 4]{ZengNewton}:
{\em At any $\tilde\bdx$ close to a semiregular zero $\bdx_*$
of a mapping $\bdx\mapsto\bdf(\bdx)$, the point $\tilde\bdx$ satisfies
the $\bdf_\bdx(\tilde\bdx)_\rkr^\dagger \bdf(\tilde\bdx) = \bdo$ 
if and only if $\bdf(\tilde\bdx)=\bdo$.}
Consequently, the {\em stationary equation} 
$\bdf_\bdx(\bdx)_\rkr^\dagger \bdf(\bdx)=\bdo$ does not produce 
extraneous zeros of $\bdf$ near any semiregular zero.

\vspace{-4mm}
\section{The low-rank Newton's iteration}

\vspace{-4mm}
It may come as a surprise that Newton's iteration we have known is only a 
special case and finding regular isolated solutions is a small portion 
of its capabilities.
The Simpson's formulation
\begin{equation}\label{nt}
\bdx_{j+1} \,=\, \bdx_j - \bdf_\bdx(\bdx_j)^{-1} \bdf(\bdx_j)
\,\,\, \mbox{for} \,\,j=0,1,\ldots
\end{equation}
is the most widely applied method for finding zeros of a mapping 
$\bdf : \Omega\subset\cV\rightarrow\cW$ if the equation $\bdf(\bdx)=\bdo$ 
is square (i.e. $\dm(\cV)=\dm(\cW)$) and the Jacobian is invertible at the 
solution.
Newton's iteration in the form of (\ref{nt}) is not suitable for computing 
singular solutions.
Even if it converges to a singular solution, the rate of convergence is 
usually slow and the attainable accuracy is poor.

A recently discovered {\em rank-$r$ Newton's iteration} \cite{ZengNewton}
\begin{equation}\label{rrnt}
\bdx_{j+1} \,=\, \bdx_j - \bdf_\bdx(\bdx_j)_\rkr^\dagger \bdf(\bdx_j)
\,\,\, \mbox{for} \,\,j=0,1,\ldots
\end{equation}
not only retains all the features of the version (\ref{nt}) but also 
expand the capability to equations of all three shapes (square, 
underdetermined and overdetermined) and to the mapping $\bdf$ 
whose Jacobian can be any rank $r$ at the solution.
Here in (\ref{rrnt}) the notation $\bdf_\bdx(\bdx_j)_\rkr^\dagger$ represents 
the Moore-Penrose inverse of the rank-$r$ projection of the Jacobian 
$\bdf_\bdx(\bdx_j)$.
The conventional Newton's iteration (\ref{nt}) and the Gauss-Newton iteration
are special cases of the rank-$r$ Newton's iteration when $r$ is the full 
column rank of the Jacobian.
This extension of Newton's method appears to be the first general purpose 
iteration for computing nonisolated solutions of the equation 
$\bdf(\bdx)=\bdo$.
The following lemma can be considered a universal convergence theorem
of Newton's iteration.

\vspace{-4mm}
\begin{lemma}[Convergence of Newton's Iteration\cite{ZengNewton}]\label{l:cnt}
Let \,$\bdf$ be a mapping twice continuously differentiable in an open domain
with a rank $r$ Jacobian $\bdf_\bdx(\bdx_*)$ at a semiregular zero $\bdx_*$.
For every open neighborhood ~$\Omega_1$ of $\bdx_*$, there is a neighborhood 
$\Omega_0$ of $\bdx_*$ such that, from every initial iterate 
$\bdx_0 \in \Omega_0$, the rank-$r$ Newton's iteration {\em (\ref{rrnt})} 
converges quadratically to a zero 
$\hat\bdx\in\Omega_1$ of \,\,$\bdf$ in the same branch as $\bdx_*$.
\end{lemma}

\vspace{-4mm}
Lemma~\ref{l:cnt} can be narrated in simpler terms: Assume an $m\times n$
equation $\bdf(\bdx)=\bdo$ has a $k$-dimensional solution set. 
Setting $r=n-k$, the rank-$r$ Newton's iteration (\ref{rrnt}) locally 
quadratically converges to a solution in the solution set if the solution
set is semiregular.
The geometric interpretation in \cite{ZengNewton} shows the iteration 
(\ref{rrnt}) asymptotically follows a normal line of the solution set and 
approximately converges to the solution nearest to the initial iterate 
$\bdx_0$.

In practical applications, equations are often given through empirical data
with limited accuracy.
On the other hand, singular solutions are highly sensitive and may even 
disappear when data are perturbed.
Those applications can be modeled as an equation
\begin{equation}\label{fxy0}
\bdf(\bdx,\bdy) \,=\,\bdo \,\,\,\mbox{for}\,\,\, \bdx\in\Omega
\end{equation}
at a fixed parameter value $\bdy$ representing the data where 
$(\bdx,\bdy)\mapsto\bdf(\bdx,\bdy)$ is a smooth mapping defined on a certain 
domain.
Assume the equation (\ref{fxy0}) has a semiregular solution $\bdx=\bdx_*$
at a data point $\bdy=\bdy_*$ but $\bdy_*$ is known only through empirical 
data $\tilde\bdy\approx\bdy_*$.
We can compute a semiregular zero of the mapping 
$\bdx\mapsto\bdf(\bdx,\bdy_*)$ near $\bdx_*$ through the perturbed 
rank-$r$ Newton's iteration 
\begin{equation}\label{niy}
\bdx_{k+1} \,=\, \bdx_k - 
\bdf_\bdx(\bdx_k,\,\tilde\bdy)_\rkr^\dagger\,\bdf(\bdx_k,\,\tilde\bdy), 
~~~~k\,=\,0,1,\ldots.
\end{equation}
If it converges, the iteration (\ref{niy}) approaches a stationary point 
$\tilde\bdx$ where 
\[
\bdf_\bdx(\tilde\bdx,\tilde\bdy)_\rkr^\dagger 
\bdf(\tilde\bdx,\tilde\bdy) = \bdo.
\]
but generally $\bdf(\tilde\bdx,\tilde\bdy) \ne \bdo$.
The following lemma ensures that the {\em stationary point} approximates
an exact solution $\check\bdx$ of the equation (\ref{fxy0}) at the exact 
data $\bdy=\bdy_*$.

\vspace{-4mm}
\begin{lemma}[Convergence of Newton's iteration on Perturbed Data 
\cite{ZengNewton}]\label{l:cpe}
Let a mapping $(\bdx,\,\bdy)\,\mapsto\,\bdf(\bdx,\,\bdy)$ be twice 
continuously differentiable in an open domain.
Assume $\bdx_*$ is a semiregular zero of the mapping 
$\bdx\,\mapsto\,\bdf(\bdx,\bdy_*)$ 
at a fixed $\bdy_*$ with $\rank{\bdf_\bdx(\bdx_*,\bdy_*)}=r>0$ ~and
$\|\bdf_\bdy(\bdx_*,\bdy_*)\| > 0$.
Then there exist a neighborhood $\Omega_*\times\Sigma_*$ of $(\bdx_*,\bdy_*)$, 
a neighborhood $\Omega_0$ of $\bdx_*$ and a constant $h$ with $0<h<1$ such 
that, at every fixed $\tilde\bdy\in\Sigma_*$ serving as empirical data for
$\bdy_*$ and from any initial iterate $\bdx_0\in\Omega_0$, the iteration 
{\em (\ref{niy})} converges to a stationary point $\tilde\bdx\in\Omega_*$ 
at which $\bdf_\bdx(\tilde\bdx,\tilde\bdy)_\rkr^\dagger\,
\bdf(\tilde\bdx,\tilde\bdy) = \bdo$ with an error bound
\begin{align}\label{eby}
\|\tilde\bdx-\hat\bdx\| ~~\le~~   \mbox{$\frac{8}{1-h}$}\,
\big\|\bdf_\bdx(\bdx_*,\bdy_*)^\dagger\big\|\,
\big\|\bdf_\bdy(\bdx_*,\bdy_*)\big\|\,\|\tilde\bdy-\bdy_*\|
+O\big(\|\tilde\bdy-\bdy_*\|^2\big) 
\end{align}
to a semiregular zero $\hat\bdx$ of $\bdx\mapsto\bdf(\bdx,\bdy_*)$ in the 
same branch of $\bdx_*$.
The convergence rate is quadratic if $\tilde\bdy=\bdy_*$.
\end{lemma}

\vspace{-4mm}
In other words, the rank-$r$ Newton's iteration (\ref{niy}) is a de facto 
{\em regularization method} that solves the exact equation (\ref{fxy0}) at 
$\bdy=\bdy_*$ approximately from perturbed data $\bdy=\tilde\bdy$.
Even though the solution of the system $\bdf(\bdx,\tilde\bdy)=0$
is substantially altered by the data perturbation or disappears altogether, 
the iteration (\ref{niy}) still converges to a stationary point $\tilde\bdx$
satisfying
$\bdf_\bdx(\tilde\bdx,\tilde\bdy)_\rkr^\dagger \bdf(\tilde\bdx,\tilde\bdy)
= \bdo$
and $\tilde\bdx$ is an accurate solution to the underlying equation
$\bdf(\bdx,\bdy_*) = \bdo$ we intend to solve. 
The accuracy of the approximate solution $\tilde\bdx$ is guaranteed by the
error bound (\ref{eby}) that is asymptotically proportional to the data error.
Furthermore, the error bound (\ref{eby}) leads to a sensitivity 
\begin{align}
\left\{\begin{array}{cl}
\big\|\bdf_\bdx(\bdx_*,\bdy_*)^\dagger\big\|_2\,\|\bdf_\bdy(\bdx_*,\bdy_*)\|_2
& \mbox{if ~$\bdx_*$ ~is semiregular} \\
\infty & \mbox{otherwise.} \end{array}\right.\label{cdn}
\end{align}
that serves as the condition number of the singular solution $\bdx_*$.
As a result, the singular zero-finding problem for 
$\bdx\mapsto\bdf(\bdx,\bdy_*)$ at empirical data $\tilde\bdy$ is regularized
as a well-posed problem of finding a stationary point $\tilde\bdx$
with an accuracy in the same order as the data.

{\bf Remark} ({\sc On identifying the projection rank})
Applying the iterations (\ref{rrnt}) and (\ref{niy}) requires identifying
the rank of the Jacobian at the zero without knowing the exact location of 
the zero or exact data of the problem.
We shall show in case studies that this rank can be determined analytically
as a part of the modeling process.
Finding the rank of a matrix from empirical data is a subject in numerical 
linear algebra as the rank-revealing problem (see, e.g. 
\cite{utvtool,li-zeng-03,lee-li-zeng,LeeLiZeng09}). 

\vspace{-4mm}
\section{Singular linear equations}

\vspace{-4mm}
Solving linear systems in the matrix-vector form $A\,\bdx=\bdb$ is one of the
most fundamental tasks in scientific computing but singular systems are rarely 
mentioned in the literature beyond elementary linear algebra.
That is an entire class of linear equations missing in discussion.
The textbook advise \cite[pp 217-218]{Meyer} is to ``avoid floating-point
solutions of singular systems'' altogether because it is well-known that 
the system becomes nominally nonsingular but highly ill-conditioned
under infinitesimal data perturbations.
As a result, the convention is to define the condition number as infinity
for singular linear systems.
In reality, however, the hypersensitivity of singular linear system is 
``notable for exaggerated fears'' \cite{PetWil79}.

From exact data, the solution of a singular linear equation $A\,\bdx=\bdb$ is 
known to be either the empty set or an affine subspace 
\begin{equation}\label{AbK}
A^\dagger \bdb + \cK(A) \,:=\,\big\{ A^\dagger \bdb + \bdz
\,\big|\, A\,\bdz=\bdo\big\}
\end{equation}
assuming $\bdb\in\cR(A)$.
Solving singular linear systems in exact sense is an ill-posed problem since
the solution generically dissipates to an empty set under arbitrary 
perturbations.
On the other hand, every vector in the affine subspace is a semiregular zero
of the mapping $\bdx\mapsto A\,\bdx-\bdb$ since the dimension of the affine 
subspace is identical to $\nullity{A}$.
Consider the holomorphic mapping 
\[\mapform{\bdf}{\CC^n\times\CC^{m\times n}\times\CC^m}{\CC^n}{(\bdx,G,\bdz)}{
G\,\bdx-\bdz.}
\]
At exact data $G=A$ of rank $r$ and $\bdz=\bdb$, the zeros of 
$\bdx\mapsto\bdf(\bdx,A,\bdb)$ form the affine subspace (\ref{AbK}).
If exact $A$ and $\bdb$ are unknown but given through empirical data
$G = \tilde{A}$ and $\bdz=\tilde{\bdb}$, the one-step rank-$r$ Newton's
iteration (\ref{niy}) from a initial iterate $\bdx_0$ produces 
\begin{align}  \tilde\bdx \,=\,
\bdx_0 + \tilde{A}_\rkr^\dagger \big(\tilde{A}\bdx_0-\tilde\bdb\big)
\,=\, \tilde{A}_\rkr^\dagger \tilde\bdb + 
\big(I-\tilde{A}_\rkr^\dagger \tilde{A}\big)\,\bdx_0 \label{tldx}
\end{align}
that accurately approximates the exact solution
\begin{equation}\label{hatx}
\hat\bdx\,=\, A^\dagger\bdb + \big(I-A^\dagger A\big)\,\bdx_0
\end{equation}
of the underlying equation $A\,\bdx = \bdb$ we intend to solve.
Moreover, the particular exact solution $\hat\bdx$ is the nearest point
in the affine subspace (\ref{AbK}) to the initial iterate $\bdx_0$.
The condition number (\ref{cdn}) is a moderate multiple of 
$\big\|A^\dagger\|_2$ that is finitely bounded and can even be small in 
applications.
The ``fear'' of singularity is indeed ``exaggerated''.

In a recent paper \cite{Zeng2019}, this author elaborates the sensitivity of
singular linear systems from a different perspective: 
The general solution of a singular system $A\,\bdx=\bdb$ is a unique point 
in an affine Grassmannian in which the sensitivity 
$\|A\|_2\big\|A^\dagger\big\|_2$ is bounded.
A properly formulated approximage solution from empirical data within
an error tolerance uniquely exists in the same affine Grassmannian, 
enjoys Lipschitz continuity and accurately approximates the exact solution
with an accuracy in the same order of the data. 
Furthermore, one can solve the perturbed system $\tilde{A}\,\bdx=\tilde\bdb$
using any method as long as it is backward stable.
The resulting solution accurately approximates one of the infinitely
many (vector) solutions.
The perceived ``errors'' are actually a part of the solution and not error 
at all.
Those results are summarized below.

\begin{theorem}[Regularization of Singular Linear Systems]
Let $A\in\CC^{m\times n}$ of rank $r$ and $\bdb\in\cR(A)$.
Assume the empirical data $(\tilde{A},\tilde\bdb)$ of $(A,\bdb)$ is
accurate so that
$\big\|\tilde{A}-A\big\|_2 < 0.46\|A^\dagger\|_2^{-1}$.
Then the following assertions hold.

{\em (i) \cite[Theorem 8]{Zeng2019}} 
Any backward accurate solution $\check\bdx$ of the data system
$\tilde{A}\,\bdx=\tilde\bdb$ is an accurate approximation to a solution
$\bdx_*$ of the underlying system $A\,\bdx=\bdb$ with an error bound
\begin{equation}\label{lserr1}
\frac{\|\check\bdx-\bdx_*\|_2}{\|\bdx_*\|_2}  \,\le\,
\mbox{\small $
\frac{\|A\|_2\,\big\|A^\dagger\big\|_2}{1-\big\|A^\dagger\|_2\,\|\Dl A \|_2}$}
\left( 2\,\sqrt{2}\,\frac{\|\Dl A\|_2}{\|A\|_2}
+\frac{\|\Dl\bdb+\bde\|_2}{\|\bdb\|_2}\right)
\end{equation}
where $\Dl A = \tilde{A}-A$, $\Dl \bdb = \tilde\bdb-\bdb$ 
and $\bde = \tilde{A}\,\check\bdx-\tilde\bdb$.

{\em (ii)} 
For any given $\bdx_0\in\CC^n$, the vector $\tilde\bdx$ in {\em (\ref{tldx})} 
is an accurate approximation to $\hat\bdx$ in {\em (\ref{hatx})} that is the 
solution of $A\,\bdx=\bdb$ nearest to $\bdx_0$ with an error bound
\begin{equation}\label{lserr2}
\big\|\tilde\bdx-\hat\bdx\big\|_2 \,\le\, \al \big\|A^\dagger\big\|_2
\big(\big\|\tilde{A}-A\big\|_2+\big\|\tilde\bdb-\bdb\big\|_2\big)+h.o.t.
\end{equation}
where $\al>0$ is a constant of moderate magnitude and $h.o.t.$ represents
higher order terms of data error.

{\em (iii) \cite[Corollary 7]{Zeng2019}} 
The affine subspace $\tilde{A}_\rkr^\dagger\tilde\bdb+
\cK\big(\tilde{A}_\rkr\big)$ accurately approximates the general solution 
{\em (\ref{AbK})} with an error bound
\begin{align}
&\max\big\{\big\|\tilde{A}_\rkr^\dagger\tilde\bdb-A^\dagger\bdb\big\|_2,\,
\dist{\cK\big(\tilde{A}_\rkr\big),\,\cK(A)} \big\} \nonumber \\
& \le\, \mbox{$
\|A\|_2\,\|A^\dagger\|_2 \frac{\sqrt{4\,\|A^\dagger\,\bdb\|_2^2+1}}{
\|A\|_2 -\|A\|_2\,\|A^\dagger\|_2\,\|A-\tilde{A}\|_2}$}\,
\big\|(\tilde{A},\,\tilde\bdb)-(A,\,\bdb) \big\|
\end{align}
\end{theorem}

For a comprehensive discussion on solving singular linear systems from 
empirical data and error analyses, see \cite{Zeng2019}.

\vspace{-4mm}
\section{Software implementation}

\vspace{-4mm}
The methods in this paper are implemented in the software package {\sc NAClab} 
\cite{naclab} on the MATLAB platform along with an intuitive interface 
\cite{solve} for solving linear and nonlinear system of equations directly as 
zero-finding for mappings, bypassing the process of representing the system 
in multivariate or matrix forms.

Solving general linear equation $L(\bdx) = \bdb$ for any linear mapping $L$, 
including singular and homogeneous cases, is implemented as the module 
{\tt LinearSolve} with an optional input item to set the error tolerance 
$\theta>0$ so that the module attempts to solve the linear equation 
$L_\theta(\bdx)=\bdb_\theta$ where $L_\theta$ is the mapping the nearest to 
$L$ with the smallest rank of all mappings within $\theta$ of $L$ and 
$\bdb_\theta$ is the orthogonal projection of $\bdb$ on $\cR(L_\theta)$.

The general Newton's iteration including the Gauss-Newton and low-rank 
Newton's iterations for solving equations in the form of $\bdf(\bdx)=\bdo$ 
is implemented as the module {\tt Newton} with the projection rank $r$ of 
the Jacobian as an optional input.

Both modules accept mappings $L$ and $\bdf$ directly as in-line MATLAB 
anonymous functions with no need to write subroutines in most cases.
Matrix representations for $L$ and the Jacobian $\bdf_\bdx(\bdx_0)$ are 
generated automatically as internal process so users can avoid the tedious 
and error-prone tasks of constructing such matrices.
We shall present several computing demos in this paper.

\vspace{-4mm}
\section{Application: Numerical Algebraic Geometry with empirical data}

\vspace{-4mm}
Numerical algebraic geometry and its application in kinematics heavily 
involve computing solutions of positive dimensions of
polynomial systems \cite{BSHW13,som-wam,WamSom11}.
Mechanisms have been developed in solving those systems including adding 
auxiliary equations to isolate witness points on the solution sets. 

When the system is given through empirical data, however, the nonisolated
solutions generally dissipates into isolated points.
The low-rank Newton's iteration can serve as a regularization mechanism
and recover the lost solutions of the underlying system.

\vspace{-2mm}
\begin{example}[Recovering lost solutions of positive dimensions]\label{e:nag}
\em
Consider the \linebreak
given mapping ~$\tilde\bdf \,:\, \CC^3\,\longrightarrow\,\CC^3$ \,defined as
\begin{align*}
\tilde\bdf ~:~  (x,&y,z) ~\longmapsto ~\\
&{\scriptsize \begin{array}{l}
(4.899\, x^3 y -5.6568\, x^5 - 8.4852\, x^3 y^2 - 2.8284\, x^3 z^2
+ 4\, x^4 + 6\, x^2 y^2 + 2\, x^2 z^2 + 7.3485\, x y^3 - 2\, x^2 \\ 
+ 2.4495\, x y z^2 + 2.8284\, x^3 - 3.4642\, x^2 y - 5.1963\, y^3 
- 1.7321\, y z^2 - 2.4495\, x y + 1.7321\, y, \\
\,\,\,\,\,\,\,\,\,\,\,\,\,\,\,\,
8.4852\, x^3 y^2 -9.798\, x^5 y - 14.697\, x^3 y^3 - 4.899\, x^3 y z^2 
+ 5.6568\, x^5 + 2.8284\, x^3 z^2 + 4.899\, x^3 y 
\\
+ 3.4642\, x^2 y z + 5.1963\, y^3 z + 1.7321\, y z^3 - 2.8284\, x^3 - 2\, x^2 z 
- z^3 - 1.7321\, y z + z - 3\, y^2 z , \\
\,\,\,\,\,\,\,\,\,\,\,\,\,\,\,\,
5.6568\, x^5 z^3 - 5.6568\, x^5 z + 5.1963\, y^3 z^2 + 1.7321 y\, z^4 
+ 16.97 x^5 y^2 z + 3.4642 x^2 y z^2+ 5.6568 x^5 
\\
+2\, x^2 z^2 -11.314\, x^7 - 4\, x^4 z^2 - 2\, x^2 z^4 + 2\, x^2 z^3 
+ 11.314\, x^7 z 
- 5.6568\, x^5 z^2 - 16.97\, x^5 y^2 - 1.7321\, y z^2 \\
- 14.697\, x^3 y^3 z - 4.899\, x^3 y z^3 + 4.899\, x^3 y z - 9.798\, x^5 y z 
- 6\, x^2 y^2 z^2 + 6\, x^2 y^2 z - 2\, x^2 z - 4.899\, x^3 y 
\\
+ 9.798\, x^5 y + 14.697\, x^3 y^3 - 5.1963\, y^3 z - 1.7321\, y z^3 
+ 1.7321\, y z 
+ 4.899\, x^3 y z^2 - 3.4642\, x^2 y z +4\, x^4 z )  \nonumber
\end{array}}
\end{align*}
as empirical data for the equation $\bdf(x,y,z)=(0,0,0)$ that is a variation of 
the system given in \cite[p. 143]{BSHW13} by replacing $x$ and $y$ with 
$\sqrt{2} x$ and $\sqrt{3} y$. 
The underlying mapping $\bdf$ has zero sets
\[ \mbox{
{\scriptsize $\big\{\sqrt{3} y = 2 x^2, z=2\sqrt{2} x^3\big\}$}, 
{\scriptsize $\big\{\sqrt{2} x = \pm 1,\sqrt{3} y = 1\big\}$},
{\scriptsize $\big\{\sqrt{2} x = 1,z = 1\big\}$}, and
{\scriptsize $\big\{2x^2+3y^2+z^2 = 1\big\}$}}
\]
that are semiregular except at intersection points.
We experiment solving this system from the data mapping $\tilde\bdf$
obtained by rounding the polynomial coefficients in five digits.

The solutions of dimension 1 and 2 disappear from rounding errors in 
coefficients. 
An attempt to solve the data system directly by Maple using rational 
coefficients did not receive results in several hours. 
The homotopy method (implemented in {\sc NAClab} as the module {\tt psolve})
terminates in seconds but results in different number 
(35-43) of isolated solutions with many of them ill-conditioned.
Bertini \cite{BSHW13} produced 36 regular solutions using hardware 
precision and 77 with adaptive multiple precision. 
The solution varieties of dimension 1 and 2 are lost by data perturbation
even if the precision is extended in floating point arithmetic.
Accurate recovery of those solutions from the given system becomes 
the problem of solving (underlying) singular system from empirical data.

By setting the projection rank $r=1$ or $r=2$, the rank-$r$ Newton's iteration
on $\tilde\bdf$ locally converges to solutions of dimension $3-r=2$ or $1$ 
respectively.
For instance, we proactively seek a solution of dimension 2 by setting $r=1$ 
from a random initial iterate.
The following is a {\sc NAClab} demo of the MATLAB command-line process 
that is intuitive without the need to write a single subroutine.

\vspace{1mm}
\scriptsize
\noindent 
\!\!$~~~${\verb|>> P = {'4.899*x^3 y - 8.4852*x^3*y^2,...;   |}~~\blue{\tt \% enter polyn. as character strings}
\newline $~~~${\verb|>> v = {'x';'y';'z'};                       |}~~\blue{\tt \% enter cell array of variable names}
\newline $~~~${\verb|>> J = PolynomialJacobian(P,v);   |}~~\blue{\tt \% Jacobian of P w.r.t. the variable names in v}
\newline $~~~${\verb|>> f = @(x,P,J,v) PolynomialEvaluate(P,v,x); |}~~\blue{\tt \% function handle 4 evaluate P at v}
\newline $~~~${\verb|>> fjac = @(x,x0,P,J,v) PolynomialEvaluate(J,v,x0)*x;  |}~~\blue{\tt \% func. evaluating J at v}
\newline $~~~${\verb|>> domain = ones(3,1); param = {P,J,v};    |}~~\blue{\tt \% domain (3x1 vectors) and parameters}
\newline $~~~${\verb|>> z0 = [-0.25518; -0.60376; -0.020624];                |}~~\blue{\tt \% random initial iterate}
\newline $~~~${\verb|>> [z,res,fcond] =  Newton({f,domain,param},{fjac,1},z0,1);      |}~~\blue{\tt \% rank-1 Newton} 
\newline $~~~~${\verb|                                     |}~~\blue{\tt \%   iteration from z0 using display type 1}

\vspace{-1.5mm}
\begin{verbatim}
Step   0:  residual =   3.59e-01 
Step   1:  residual =   4.67e-02    shift =   4.99e-02
Step   2:  residual =   1.25e-03    shift =   8.88e-03
Step   3:  residual =   9.74e-07    shift =   2.51e-04
Step   4:  residual =   6.93e-08    shift =   1.96e-07
Step   5:  residual =   6.93e-08    shift =   1.20e-13
Step   6:  residual =   6.93e-08    shift =   8.50e-17
\end{verbatim}
\normalsize

\vspace{1mm}
Notice that the residual can only reduce to $6.93\times 10^{-8}$. 
Namely the limit $(\tilde{x},\tilde{y},\tilde{z})$ is {\em not} a zero of 
$\tilde\bdf$ but a stationary point as a solution to
$\tilde\bdf_{xyz}(\tilde{x},\tilde{y},\tilde{z}
)_\rk{1}^\dagger\bdf(\tilde{x},\tilde{y},\tilde{z}) = \bdo.$
as indicated by the shifts 
\[
\|(x_{j+1},y_{j+1},z_{j+1})-(x_j,y_j,z_j)\|_2, \,\,\,
j=0,1,\ldots
\]
approaching hardware zero.
The stationary equation regularizes the singular equation $\bdf(x,y,z) = \bdo$.
The iteration terminates in 6 steps at 
\[
(\tilde{x},\tilde{y},\tilde{z}) = (\mbox{\scriptsize
$-0.234036969240715  -0.544684891672585  -0.020211408075956$})
\]
that accurately approximates a point $(\check{x},\check{y},\check{z})$
in the solution set $\{2x^2+3y^2+z^2=1\}$ with 10 correct digits. 
\end{example}

\vspace{-4mm}
\section{Application: the GCD equation}

\vspace{-4mm}
An intuitive model for computing the greatest common divisor (GCD) of a 
polynomial pair $p$ and $q$ is solving the GCD equation
\begin{equation}\label{gcdeq}
  \big(u\,v - p, \,\,u\,w-q\big) \,\,=\,\, (0,\,0)
\end{equation}
for $(u,v,w)=(u_*,v_*,w_*)$ where $u_*$ is a constant multiple of the GCD
and $(v_*,w_*)$ is a pair of co-factors.
However, the equation (\ref{gcdeq}) is obviously singular with the 
1-dimensional solution set 
\begin{equation}\label{gcdsol}
\big\{ 
\big(t\,u_*,\,\mbox{$\frac{1}{t}$}\,v_*,\,\mbox{$\frac{1}{t}$}\,w_*
\big)
\,\big|\, t\in\CC\setminus\{0\}\big\}.
\end{equation}
The model (\ref{gcdeq}) is viable only if we can accurately solve for its
singular solutions, or after adding unnatural auxiliary constraints.
Furthermore, the solution set (\ref{gcdsol}) is infinitely sensitive and 
generically reduces to the trivial set 
$\big\{\big(t,p/t,q/t)\,\big|\,t\in\CC\setminus\{0\}\big\}$
under arbitrary data perturbations.
As a result, solving the equation (\ref{gcdeq}) has been an ill-posed problem
with empirical data.
On the other hand, the emergence of the low-rank Newton's iteration enables 
solving the GCD equation (\ref{gcdeq}) directly and accurately even if the
data are perturbed and the nontrivial GCD disappears.

Let $\cP_l$ denote the vector space of polynomials with degrees up to $l$ with
a norm $\|u\|$ defined as the 2-norm of the coefficient vector of $u\in\cP_l$
so that $\cP_l$ is isometrically isomorphic to $\CC^{l+1}$.
Assume $(p,q)\in\cP_m\times\cP_n$ of degrees $m$ and $n$, respectively, 
with the GCD degree $k$. 
We define the holomorphic mapping 
\begin{equation}\label{gcdmap}
\mapform{\bdf}{\cP_k\times\cP_{m-k}\times\cP_{n-k}\times\cP_m\times\cP_n}{
\cP_m\times\cP_n}{(u,v,w,f,g)}{\big(u\,v-f,\,u\,w-g\big)}
\end{equation}
The following lemma establishes the semiregularity of the solution set 
(\ref{gcdsol}).

\begin{lemma}[Semiregularity of the GCD Equation]\label{l:gcd}
Let $\bdf$ be defined in {\em (\ref{gcdmap})}.
Assume $p$ and $q$ are polynomials of degrees $m$ and $n$, respectively,
with the GCD degree $k$.
Then any zero $\big(\hat{u},\hat{v},\hat{w}\big)$ of the mapping 
$\bdg : (u,v,w)\mapsto\bdf(u,v,w,p,q)$ at the fixed parameter value 
$(f,g) = (p,q)$ with $\dg(\hat{u})=k$ is semiregular.
\end{lemma}

{\em Proof.} Let $(u_*,v_*,w_*)$ be a particular zero of $\bdg$ with
$\dg(u_*)=k$. 
Then (\ref{gcdsol}) is the zero set $\bdg^{-1}(\bdo)$ of dimension one.
Consequently $\nullity{\bdg_{uvw}(\hat{u},\hat{v},\hat{w})}\ge 1$.
By \cite[Lemma 4.1]{ZengGCD}, appending one extra linear equation to
$\bdg(u,v,w)=\bdo$ makes the Jacobian of the left side injective at 
$(\hat{u},\hat{v},\hat{w})$. 
Namely $\nullity{\bdg_{uvw}(\hat{u},\hat{v},\hat{w})} \le 1$.
As a result, the zero $(\hat{u},\hat{v},\hat{w})$ is semiregular since
the nullity of $\bdg_{uvw}(\hat{u},\hat{v},\hat{w})$ is 1 and identical
to the dimension of (\ref{gcdsol}) \qed

Since the zero set (\ref{gcdsol}) of $\bdg$ is of dimension 1 and semiregular,
the rank of the Jacobian at any particular solution is 
\begin{equation}\label{gcdrank}
 r\,=\,\dm\big(\cP_k\times\cP_{m-k}\times\cP_{n-k}\big)-1 \,=\,m+n-k+2 
\end{equation}
and the rank-$r$ Newton's iteration with $r$ as in (\ref{gcdrank}) 
\begin{align}
(u_{j+1},v_{j+1},w_{j+1})\,=\,(u_j,v_j,w_j) 
&-\bdf_{uvw}(u_j,v_j,w_j,\tilde{p},\tilde{q})_\rkr^\dagger\,
\bdf(u_j,v_j,w_j,\tilde{p},\tilde{q})\label{gcdit}
\end{align}
for $j=0,1,\ldots$ at empirical data $(f,g)=(\tilde{p},\tilde{q})$
locally converges to a GCD triple $(\tilde{u},\tilde{v},\tilde{w})$ that 
accurately approximates an exact GCD triple $(\hat{u},\hat{v},\hat{w})$ of 
the underlying data $(p,q)$ in the zero set (\ref{gcdsol}).
Consequently, the iteration (\ref{gcdit}) serves as an effective regularization
mechanism for the singular GCD equation (\ref{gcdeq}) so that the GCD can
be accurately computed from empirical data, 
as asserted in the following theorem.

\begin{theorem}[Regularization of GCD] 
Let $(p,q)$ be a polynomial pair of degrees $m$ and $n$ respectively with
a GCD degree $k$.
Assume the data $(\tilde{p},\tilde{q})$ is sufficiently close to $(p,q)$
and the initial iterate $(u_0,v_0,w_0)$ is sufficiently close to a zero
$(\hat{u},\hat{v},\hat{w})$ of the mapping
$\bdg : (u,v,w)\mapsto\bdf(u,v,w,p,q)$.
Setting $r=m+n-k+2$, the rank-$r$ Newton's iteration {\em (\ref{gcdit})}
at the data $(\tilde{p},\tilde{q})$ converges to 
$(\tilde{u},\tilde{v},\tilde{w})$ with an error bound
\begin{align}
\big\|(\tilde{u},\,&\tilde{v},\tilde{w})-(\check{u},\check{v},\check{w})\big\|
\,\le 
&c\,\big\|\bdf_{uvw}(\hat{u},\hat{v},\hat{w},p,q)_\rkr^\dagger\big\|
\,\|(\tilde{p},\tilde{q})-(p,q)\| + h.o.t. \label{gcderr}
\end{align}
where $(\check{u},\check{v},\check{w})\in\bdg^{-1}(\bdo)$ is an exact
GCD triple of $(p,q)$ and $c>0$ is a constant of moderate magnitude.
The convergence rate is quadratic if data $(\tilde{p},\tilde{q})=(p,q)$.
The GCD condition number can be defined as
\,$\|\bdf_{uvw}(\hat{u},\hat{v},\hat{w},p,q)_\rkr^\dagger\|$
\,at the polynomial pair $(p,q)$ 
\end{theorem}

{\em Proof.} A straightforward verification based on 
Lemma~\ref{l:cpe} and Lemma~\ref{l:gcd}. \qed

The GCD model (\ref{gcdeq}) is not restricted to the univariate GCD problem.
Multivariate GCD's can be computed by solving the same equation from proper 
domains of polynomial spaces using the same iteration (\ref{gcdit}) except
that the projection rank $r$ needs to be adjusted to one less than the 
dimension of the corresponding domain.

\vspace{-4mm}
\section{Application: Factoring polynomials}

\vspace{-4mm}
A straightforward and intuitive model for factoring a multivariate polynomial 
$p$ is to solve the factorization equation 
\begin{equation}\label{faceq}
u_0\,u_1^{\ell_1}\,\cdots\,u_k^{\ell_k} - p \,=\, 0
\end{equation}
for an irreducible factor array $(u_0,\ldots,u_k)=(\hat{u}_0,\ldots,\hat{u}_k)$
where $\ell_1,\ldots,\ell_k\,>\,0$ are integers.
For convenience, we assume $\hat{u}_0\in\CC$, $\ell_0=1$ and 
$\hat{u}_1,\ldots,\hat{u}_k$ are nontrivial.
The equation (\ref{faceq}) is singular with a solution set of dimension 
$k$ in the form of
\begin{equation}\label{facsol}
\big\{\big(t_0 \hat{u}_0,\,t_1 \hat{u}_1,\ldots,\, t_k \hat{u}_k\big)
\,\big|\, t_1,\ldots,t_k\in\CC\!\setminus\!\{0\}, 
t_0 = t_1^{-\ell_1}\cdots t_k^{-\ell_k}
\big\}
\end{equation}
which is hypersensitive and the exact nontrivial factorization is generally
impossible if $p$ is known only through empirical data $\tilde{p}$.

Let $\cU_0=\CC$ and $\cU_1$, $\ldots$, $\cU_k$ be vector spaces of polynomials 
containing $\hat{u}_0,\ldots,\hat{u}_k$ respectively.
For $j=0,1,\ldots,k$, assume every $\cU_j$ is a {\em proper hosting space} 
of $\hat{u}_j$ in the sense that $s\,\hat{u}_j\in\cU_j$ implies $s$ is a 
constant.
Let $\cP$ be a vector space of polynomials containing $p$, $\tilde{p}$ and 
all the products $u_0 u_1\cdots u_k$ for $u_j\in\cU_j$, $j=0,1,\ldots,k$.
Define the holomorphic mapping
\begin{equation}\label{facmap}
\mapform{\bdf}{\cU_0\times\cU_1\times\cdots\times\cU_k\times\cP}{\cP}{
(u_0,u_1,\ldots,u_k,f)}{u_0\, u_1^{\ell_1}\cdots u_k^{\ell_k}-f}
\end{equation}
The following lemma establishes the crucial semiregularity of 
(\ref{facsol}).

\begin{lemma}[Semiregularity of Polynomial Factorization]\label{l:fac}
Let \,$\bdf$ be defined in {\em (\ref{facmap})} and $p\in\cP$ with an
irreducible factorization $\hat{u}_0\hat{u}_1^{\ell_1}\cdots\hat{u_k}^{\ell_k}$ 
where $\hat{u}_j$ belongs to a proper hosting space 
$\cU_j$ for $j=0,1,\ldots,k$.
Then every zero $(\check{u}_0,\ldots,\check{u}_k)$ of the mapping
$\bdg : (u_0,\ldots,u_k)\mapsto\bdf(u_0,\ldots,u_k,p)$ at 
$f=p$ is semiregular and
\begin{equation}\label{facr}
  \rank{\bdg_{u_0\cdots u_k}(\check{u}_0,\ldots,\check{u}_k)}
\,=\,\dm\big(\cU_0\times\cdots\times\cU_k\big)-k
\end{equation}
\end{lemma}

{\em Proof.} The Jacobian 
$\bdg_{u_0\cdots u_k}(\check{u}_0,\ldots,\check{u}_k)$ is the linear map
\[
(u_0,\ldots,u_k) \mapsto \mbox{
$\sum_{i=0}^k u_i \big(\ell_i \check{u}_i^{\ell_i-1} 
\prod_{j\ne i} \check{u}_j^{\ell_j}\big)$}
\]
whose nullity is at least $k$ since the zero set 
(\ref{facsol}) of $\bdg$ is of dimension $k$.
Let $\phi_j : \cU_j\rightarrow \CC$ be a linear functional
with $\phi_j(\hat{u}_j)=\bt_j\ne 0$ for $j=1,\ldots,k$.
Consider the mapping 
\begin{align}\label{h}
\bdh ~:~ (u_0,\ldots,u_k) ~\longmapsto~
\big(\bdf(u_0,\ldots,u_k,p),\,\,\phi_1(u_1)-\bt_1,\,\,\ldots,\,\,
\phi_k(u_k)-\bt_k\big)
\end{align}
and we claim its Jacobian at $(\check{u}_0,\ldots,\check{u}_k)$ is injective.
In fact, setting 
\[
\bdh_{u_0\cdots u_k}(\check{u}_0,\ldots,\check{u}_k)(u_0,\ldots,u_k) = \bdo
\]
yields, for any $i\in\{0,\ldots,k\}$, 
\[
u_i \Big(\ell_i\,
\mbox{$\prod_{j\ne i}$} \check{u}_j^{\ell_j}\Big) \, 
=  \, -\check{u}_i \Big(\mbox{$\sum_{l\ne i} u_l \big(
\ell_l \check{u}_l^{\ell_l-1} 
\prod_{j\ne l,i}$}\check{u}_j^{\ell_j}\big) \Big)
\]
implying $u_i=s\check{u}_i$ and $s$ must be a constant.
As a result, we have $\phi_i(s\check{u})=s\phi_i(\check{u})=0$, leading to 
$s=0$.
Thus $u_i=0$ for all $i=0,\ldots,k$ so the Jacobian of $\bdh$ is injective 
at $(\check{u}_0,\ldots,\check{u}_k)$.
Since appending $k$ linear functionals to 
$\bdg_{u_0\cdots u_k}(\check{u}_0,\ldots,\check{u}_k)$
reduces its nullity to zero, its nullity is no more than $k$, leading to the
semiregularity of $(\check{u}_0,\ldots,\check{u}_k)$ and (\ref{facr}) holds. 
\qed

Setting $r$ as (\ref{facr}) by Lemma~\ref{l:fac}, the rank-$r$ Newton's 
iteration 
\begin{align}\label{facit}
\big(u_0^{(j+1)},\ldots,u_k^{(j+1)}\big) &\,=\, 
\big(u_0^{(j)},\ldots,u_k^{(j)}\big)   \\ &
-\bdf_{u_0\cdots u_k}\big(
u_0^{(j)},\ldots,u_k^{(j)},\tilde{p}\big)_\rkr^\dagger\,
\bdf\big(u_0^{(j)},\ldots,u_k^{(j)},\tilde{p}\big) \nonumber
\end{align}
regularizes the factorization problem as asserted in the following theorem.

\begin{theorem}[Regularization of Polynomial Factorization]\label{p:fac}
Let $p = \hat{u}_0\,\hat{u}_1^{\ell_1}\cdots\hat{u}_k^{\ell_k}$ 
be an irreducible polynomial factorization where 
$\hat{u}_j$ belongs to a proper hosting space $\cU_j$ 
for $j=0,\ldots,k$ and $\cU_0=\CC$.
Let $\cP\ni p$ be a vector space containing all products 
$u_0\,u_1^{\ell_1}\cdots u_k^{\ell_k}$ for $u_j\in\cU_j$, $j=0,\ldots,k$
and set $r$ to be {\em (\ref{facr})}. 
Then, for any $\tilde{p}\in\cP$ sufficiently close to $p$
as empirical data and from any initial iterate 
$\big(u_0^{(0)},\ldots,u_k^{(0)}\big)\in\cU_0\times\cdots\times\cU_k$ near
$(\hat{u}_0,\ldots,\hat{u}_k)$, the rank-$r$ Newton's iteration 
{\em (\ref{facit})} converges to a 
$\big(\tilde{u}_0,\ldots,\tilde{u}_k\big)\in\cU_0\times\cdots\times\cU_k$ 
with an error bound
\begin{align}\label{facerr}
\big\|\big(\tilde{u}_0,\ldots,\tilde{u}_k\big)-&
\big(\check{u}_0,\ldots,\check{u}_k\big)\big\|_2  \\
&\le ~~\al\,\big\|\bdf_{u_0\cdots u_k}(\hat{u}_0,\ldots,\hat{u}_k,p
)_\rkr^\dagger\big\| \|\tilde{p}-p\|+O(\|\tilde{p}-p\|^2) \nonumber
\end{align}
where $\big(\check{u}_0,\ldots,\check{u}_k\big)$ an exact factor array of 
$p$ in {\em (\ref{facsol})} and $\al>0$ is a constant of moderate size.
The convergence is quadratic if $\tilde{p}=p$.
The norm
\,$\|\bdf_{u_0\cdots u_k}(\hat{u}_0,\ldots,\hat{u}_k,p)_\rkr^\dagger\|$
can be defined as the factorization condition number of $p$.
\end{theorem}

{\em Proof.} The assertions follow from a straightforward verification
using Lemma~\ref{l:cpe} and Lemma~\ref{l:fac}. \qed

Regularizing the singular factorization problem by taking advantage of
the semiregularity and the low-rank Newton's iteration (\ref{niy})
substantially improves the existing results in \cite{WuZeng} theoretically
and computationally by eliminating the unnatural auxiliary components
$\phi_j(u_j)-\bt_j$ for $j=1,\ldots,k$ in (\ref{h}) from the model.

\begin{example}[Factoring a polynomial from empirical data]
\em
The data for the \newline polynomial $p=\mbox{\scriptsize
$\big(\frac{2}{3}y^2+\frac{3}{7}x^2 z^6\big)^3 
\big(-1+\frac{5}{11} y z + x^5\big)^2$}$ is given in 
\begin{align*}
\tilde{p} &~=~ \\
&{\scriptsize \begin{array}{l}
.296296 y^9 - 0.269360 y^{10} z - 1.02640 y^9 x^5 + 0.0612182 y^{11} z^2 + 0.466545 y^{10} z x^5 + 0.888889 y^9 x^{10}\\
+ 1.14286 y^6 x^2 z^4 - 1.03896 y^7 x^2 z^5 - 3.95896 y^6 x^7 z^4 + 0.236128 y^8 x^2 z^6 + 1.79953 y^7 x^7 z^5\\
+ 3.42857 y^6 x^{12} z^4 + 1.46939 y^3 x^4 z^8 - 1.33581 y^4 x^4 z^9 - 5.09011 y^3 x^9 z^8 + 0.303593 y^5 x^4 z^{10}\\
+ 2.31369 y^4 x^9 z^9 + 4.40816 y^3 x^14 z^8 + 0.629738 x^6 z^{12} - 0.572489 x^6 z^{13} y - 2.18148 x^11 z^{12}\\
+ 0.130111 x^6 z^{14} y^2 + 0.991580 x^{11} z^{13} y + 1.88921 x^{16} z^{12}
\end{array}}
\end{align*}
From the data polynomial $\tilde{p}$, the factorization structure of $p$ 
can be identified by the methods elaborated in \cite{WuZeng} along with
initial approximation of factors.
Applying Proposition~\ref{p:fac} with $k=2$, $\ell_1=3$, $\ell_2=2$
along with fewnomial spaces $\cU_1 = \spn\{y^3,x^2 z^4\}$ and 
$\cU_2 = \spn\{1, y z,x^5\}$, we can carry out the rank-4 Newton's iteration
(\ref{facit}) in the follow computing demo of {\sc NAClab} in which
{\tt pplus, pminus, ptimes} are polynomial utilities for $+$, $-$ and $\times$.

\vspace{1mm}
\scriptsize
\noindent 
\!\!$~~~${\verb|>> P = {.296296*y^9 - 0.269360*y^10*z ...;    |}~~\blue{\tt \% enter polynomials as char. strings}
\newline $~~~${\verb|>> f = @(u,v,w,p) pminus(ptimes(u,v,v,v,w,w,p);   |}~~\blue{\tt \% function handle for mapping f}
\newline $~~~${\verb|>> fjac = @(u,v,w,u0,v0,w0,p) pplus(ptimes(u,v0,v0,v0,w0,w0),ptimes(u0,3,v0,...;|}
\newline $~~~${\verb|>> v0,v,w0,w0), ptimes(u0,v0,v0,v0,2,w0,w));   |}~~\blue{\tt \% function of the Jacobian mapping}
\newline $~~~${\verb|>> domain = {1,'y^3+x^2*z^4', '1+y*z+x^5'}; param = {p};  |}~~\blue{\tt \% domain and parameters}
\newline $~~~${\verb|>> u0=1; v0='.67*y^3+.86*x^2*z^4'; w0='-1+.45*y*z+1.73*x^5';    |}~~\blue{\tt \% initial iterate}
\newline $~~~${\verb|>> [Z,res,fcnd]=Newton({f,domain,param},{fjac,4},{u0,v0,w0},1)    |}~~\blue{\tt \% rank-4 Newton}
\newline $~~~~${\verb|      |}~~\blue{\tt \%  on f with given domain, param \& Jac. from (u0,v0,w0) using display type 1}

\vspace{-1mm}
\begin{verbatim}
   Step   0:  residual =   5.35e-02 
   Step   1:  residual =   2.26e-04    shift =   3.70e-03
   Step   2:  residual =   7.87e-06    shift =   3.53e-05
   Step   3:  residual =   7.86e-06    shift =   1.22e-09
   Step   4:  residual =   7.86e-06    shift =   6.58e-16
\end{verbatim}
\normalsize

\noindent
The process terminates at the approximate factorization
\[
0.999035 \big(.667678 y^3 + .858444 x^2 z^4\big)^3
\big(-.998210 + .453732 y z + 1.7289489 x^5\big)^2
\]
toward a point in the 2-dimensional solution (\ref{facsol}) with coefficients 
accuracy $6.8\times 10^{-6}$ that is in the same order as the data accuracy.
Again, residual does not approach zero since the mapping $\bdf$ in 
(\ref{facmap}) does {\em not} have a zero at the data $f=\tilde{p}$.
However, the shifts approaching zero implies the iteration solves the
stationary equation 
\[
\bdf_{u_0 u_1 u_2}(u_0,u_1,u_2,\tilde{p})_\rk{4}^\dagger
\bdf(u_0,u_1,u_2,\tilde{p}) = 0.
\]
that regularizes the singular equation (\ref{faceq}) with a 
near optimal condition number 4.92.
The {\sc NAClab} equation solving interface \cite{solve} makes the
entire process intuitive.
\end{example}

\vspace{-4mm}
\section{Application: Defective eigenvalues}

\vspace{-4mm}
Computing defective eigenvalues of matrices is a well-known singular problem
and a formidable challenge to achieve accurate results from empirical data.
We shall demonstrate that defective eigenvalues are semiregular
and can be regularized through the low-rank Newton's iteration,
advancing from sensitivity theory and the computational method in 
\cite{pseudoeig}.

Let $\hat\la$ be a defective eigenvalue of a matrix $A\in\CC^{n\times n}$.
We say the {\em multiplicity support} of $\hat\la$ is $m\times k$ if $\hat\la$
is of geometric multiplicty $m$ with the smallest Jordan block size $k$.
Finding such an eigenvalue can naturally modeled as solving the 
eigenequation
\begin{equation}\label{eigeq}
A\,X - \la \,X - X\,S \,=\, O
\end{equation}
for $(\la,X)\in\CC\times\CC^{n\times k}$ where $S\in\CC^{k\times k}$ satisfies
\begin{equation}\label{S}
S=[s_{ij}], \,\,\, s_{ij}=0 \,\,\mbox{for $i\le j$ and} 
\,\,\,s_{12} s_{23} \cdots s_{k-1,k} \ne 0.
\end{equation}
Any solution of (\ref{eigeq}) is a zero of the mapping 
$(\la,X)\mapsto\bdf(\la,X,A)$ where
\begin{equation}\label{eigmap}
\mapform{\bdf}{\CC\times\CC^{n\times k}\times\CC^{n\times n}}{
\CC^{n\times k}}{(\la,X,G)}{G\,X-\la\,X-X\,S}
\end{equation}

\begin{lemma}[Semiregularity of Defective Eigenvalues]\label{l:eig}
Let $\hat\la$ be an eigenvalue of $A\in\CC^{n\times n}$ with a multiplicity 
support $m\times k$.
For any fixed parameter $S$ satisfying {\em (\ref{S})}, the solution of 
{\em (\ref{eigeq})} is semiregular with dimension $m\,k$ in the form of
\begin{equation}\label{eigsol}
\big\{\big(\hat\la, \,\hat{X}\big) \,\big|\, \hat{X} = X_0+Y_0\,Z,\, 
Z\in\CC^{m\times k} \big\}
\end{equation}
where $X_0,\,Y_0 \in\CC^{n\times k}$ with $\cR(Y_0)=\cK(A-\hat\la I)$.
Furthermore, the partial Jacobian $\bdf_{_{\la X}}(\hat\la,\hat{X},A)$
at any solution is of rank
\begin{equation}\label{eigr}
r\,=\,\rank{\bdf_{_{\la X}}\big(\hat\la,\,\hat{X},\,A\big)} \,=\,
1+(n-m)\,k.
\end{equation}
\end{lemma}

{\em Proof.} Write $X = \blb \bdx_1,\ldots,\bdx_k\brb$ columnwise. 
Then the equation (\ref{eigeq}) with $\la=\hat\la$ can be expanded as
$(A-\hat\la I)\,\bdx_1 \,=\,\bdo$ along with
\[
A\,\bdx_j-\hat\la\,\bdx_j \,=\,s_{1j}\bdx_1+\cdots+s_{j-1,j}\bdx_{j-1}, \,\, 
\mbox{for}\,\,j=2,\ldots,k.
\]
by picking any specific solution $X=X_0$ satisfying the above system and any 
$Y_0\in\CC^{n\times m}$ whose columns form a basis for $\cK(A-\hat\la I)$, 
we have the solution (\ref{eigsol}) since $\nullity{A-\hat\la I}=m$.
The mapping $\phi : Z\mapsto (\hat\la, X_0+Y_0\,Z)$ is injective
and the Jacobian $\phi_{_Z} (\hat{Z}) : Z\mapsto Y_0\,Z$ is of rank $m\,k$
since $Y_0$ is of full column rank.
Define the mapping 
\[ \mapform{\bdg}{\CC\times\CC^{n\times k}}{\CC^{n\times k}}{(\la,\,X)}{
(A-\la I)\,X-X\,S}
\]
Thus the Jacobian $\bdg_{_{\la X}}$ at any solution in (\ref{eigsol}) 
is of nullity at least $m\,k$.
By \cite[Lemma~2]{pseudoeig}, appending a linear mapping 
$X\mapsto C^\h X\in\CC^{m\times k}$ with a constant matrix 
$C\in\CC^{n\times m}$ to $\bdg_{_{\la X}}(\hat\la,\hat{X})$ reduces the 
nullity to zero, implying $\nullity{\bdg_{_{\la X}}(\hat\la,\hat{X})}$ is
no more than $m\,k$.
Hence every solution in (\ref{eigsol}) is semiregular.
The rank (\ref{eigr}) follows accordingly.
\qed

Upon establishing semiregularity and setting the projection rank $r$
in (\ref{eigr}), we can now compute a defective eigenvalue from empirical data 
$\tilde{A}$ by applying the rank-$r$ Newton's iteration
\begin{align}\label{eigit}
(\la_{j+1},X_{j+1}) \,=\, (\la_j,X_j) - 
\bdf_{_{\la X}}(\la_j,X_j,\tilde{A})_\rkr^\dagger \bdf(\la_j,X_j,\tilde{A})
\end{align}

\begin{theorem}[Regularization of Defective Eigenvalues]\label{p:eig}
Let $\hat\la$ be an eigenvalue of $A\in\CC^{n\times n}$ with a multiplicity
support $m\times k$.
Then, for any $\tilde{A}\in\CC^{n\times n}$ sufficiently close to $A$ as 
empirical data and $r$ as in {\em (\ref{eigr})}, the rank-$r$ Newton's 
iteration {\em (\ref{eigit})} from any initial iterate $(\la_0,X_0)$ close
to a solution in {\em (\ref{eigsol})} converges to a point
$(\tilde\la,\tilde{X})\in\CC\times\CC^{n\times k}$ with an error bound
\begin{align}\label{eigerr}
\big\|(\tilde\la,\tilde{X})&-(\hat\la,\check{X})\big\| \,\le\, 
\al\,\big\|\bdf_{_{\la X}}(\hat\la,\check{X},A)_\rkr^\dagger\big\|\,
\big\|\check{X}\big\|_2 \|A-\tilde{A}\|_{_F}
+ O(\|A-\tilde{A}\|_{_F}^2) 
\end{align}
where $\al=O(1)$ is a constant and $(\hat\la,\check{X})$ is an exact solution
in {\em (\ref{eigsol})}.
The convergence is quadratic if $\tilde{A}=A$.
The condition number of $\hat\la$ can be defined as
$\big\|\bdf_{_{\la X}}(\hat\la,\check{X},A)_\rkr^\dagger\big\|$
\end{theorem}

{\em Proof.} The assertions directly follows Lemma~\ref{l:cpe} and 
Lemma~\ref{l:gcd} with $\big\|\bdf_{_G}(
\hat\la,\check{X},A)\big\| \le \|\check{X}\|_2$ \,since
$\bdf_{_G}(\hat\la,\check{X},A) : G \mapsto G\,\check{X}$.
The $\|\hat{X}\|_2$ component in the condition number of $\hat\la$ can be
eliminated since the $\hat{X}$ can be chosen with orthonormal columns.
\qed

The error estimate (\ref{eigerr}) can be improved by eliminating the factor
$\|\check{X}\|_2$ by a thin-QR decomposition $\tilde{X}=Q\,R$, resetting
the component $S$ as $R\,S\,R^{-1}$ and one additional step of
the iteration (\ref{eigit}) from the initial iterate 
$(\la_0,X_0)=(\tilde\la, Q)$.
The resulting $X$ component will have nearly orthonormal columns and 2-norm
approximately one.
This normalization process is much simpler than that in \cite{pseudoeig}.

\begin{example}[Defective eigenvalue from empirical data]
\em
Let $\hat\la=2$ be a 7-fold eigenvalue of $A$ with multiplicity support 
$2\times 2$ but $A$ is known through data $\tilde{A}$ below
with entry error bound $.5\times 10^{-4}$.
\[  A ~=~ \mbox{\scriptsize $
\left[\begin{array}{rrrrrrrr}
-0.1047 &  2.6711 & -7.7657 &  7.6782 & -0.1741 & -2.8614 & -1.5102 & 10.1186\\
 1.1993 &  1.3389 &  2.5196 & -2.4136 & -0.5598 &  1.1995 &  1.5892 & -3.1106\\
 1.5919 & -4.4314 & 10.3181 & -7.9651 &  0.8970 &  1.3103 &  0.2183 &-11.4464\\
 3.0877 & -4.2142 &  9.8737 & -7.5953 &  0.4991 &  3.1022 &  1.4778 &-13.1894\\
 1.3996 &  0.6824 &  0.3731 & -0.3272 &  1.2337 &  0.4494 &  0.6920 & -0.0206\\
 0.2930 & -0.4477 &  1.8217 & -2.4647 & -0.3103 &  3.4128 &  0.7911 & -2.8883\\
 0.8370 & -0.3341 &  1.7179 & -0.9933 &  0.4461 &  0.2581 &  1.8852 & -1.4502\\
-1.7541 &  0.4549 & -2.9046 &  2.8613 &  0.4126 & -1.9328 & -1.5465 &  5.5124
\end{array}\right]$}
\]
Matlab built-in function {\tt eig} produces scattered eigenvalues
\[
\mbox{\scriptsize $1.7733 \pm 0.1345i,\, 2.0341 \pm 0.2668i,\,   
2.1931 \pm 0.0454i,\, 1.9976,\,   2.0025 $}
\]
of errors at least $.0024$.
From an initial estimate $\la_0 = 1.98$, we first calculate the component
$X_0$ of the initial iterate by solving 
\[
A\,X-\la_0 X-X\,S \,=\,O
\]
for $X\in\CC^{8\times 2}$ within error tolerance $3\times 10^{-2}$ in the 
following {\sc NAClab} calling sequence:

\scriptsize
\noindent 
\!\!$~~~${\verb|>> A = [-0.1047 2.6711 -7.6782 ...; |}~~\blue{\tt \% enter data matrix}
\newline $~~~${\verb|>> S = [0 1; 0 0];                  |}~~\blue{\tt \%  matrix parameter S}
\newline $~~~${\verb|>> L = @(X,e0,G,S) G*X-e0*X-X*S;    |}~~\blue{\tt \% function handle for L : X -> G*X-e0*X-X*S}
\newline $~~~${\verb|>> [~,K]=LinearSolve({L,{ones(8,2)},{1.98,A,S}},zeros(8,2),3e-2) |}~~\blue{\tt \% solve L(X)=O}

\normalsize
\noindent
obtaining the initial iterate $(\la_0,X_0)$ where $X_0$ is a random linear
combination of the four solutions in output {\tt K} of {\tt LinearSolve}.
The rank-13 Newton's iteration is carried out as follows.

\vspace{-1mm}
\scriptsize
\noindent 
\newline $~~~${\verb|>> f = @(e,X,G,S) G*X-e*X-X*S;|}~~\blue{\tt \% function handle for mapping f:(e,X)->G*X-e*X-X*S}
\newline $~~~${\verb|>> fjac = @(e,X,e0,X0,G,S) G*X-e*X0-e0*X-X*S;|}~~\blue{\tt \% Jacobian (e,X)->G*X-e*X0-e0*X-X*S}
\newline $~~~${\verb|>> domain = {1,ones(8,2)};                           |}~~\blue{\tt \% domain of f as C x C\^{ }\{8x2\}}
\newline $~~~${\verb|>> param = {A,S};                              |}~~\blue{\tt \% parameters A and S for mapping f}
\newline $~~~${\verb|>> [Z,res,fcnd]=Newton({f,domain,param},{fjac,13},{e0,X0},1);    |}~~\blue{\tt \% rank-13 Newton}

\begin{verbatim}
   Step   0:  residual =   8.78e-02 
   Step   1:  residual =   2.13e-04    shift =   2.01e-02
   Step   2:  residual =   1.36e-05    shift =   3.95e-04
   Step   3:  residual =   1.36e-05    shift =   3.47e-09
   Step   4:  residual =   1.36e-05    shift =   1.86e-14
   Step   5:  residual =   1.36e-05    shift =   1.41e-15
\end{verbatim}
\normalsize

\noindent
obtaining an accurate defective eigenvalue $\tilde\la = 2.000072$ with an
accuracy $.7\times 10^{-4}$ in the same level of the data error.
\end{example}

\vspace{-4mm}
\section{On ultrasingular equations}

\vspace{-4mm}
We say an equation is {\em ultrasingular} if its Jacobian at a desired 
solution has a (column) rank-deficiency larger than the dimension of 
the solution.
Ultrasingularity occurs in cases such as at a zero whose dimension is 
undefined (e.g. intersection points of solution branches), isolated multiple
zeros, isolated ultrasingular zeros embedded in a semiregular zero set
and entire branch of nonisolated ultrasingular zeros.
Difficulties in computing ultrasingular zeros including slow convergence
rate of iterative methods (c.f. \cite{DecKelKel}) and, more importantly, 
barriers of low attainable accuracy \cite{victorpan97,ypma}.

Singular equations with isolated multiple zeros can be accurately solved
by the {\em depth-deflation} method \cite{DLZ,dayton-zeng}:
A singular isolated zero $\bdx_*$ of a mapping 
$\bdf : \Omega\subset\CC^m\rightarrow\CC^n$ derives an isolated zero
$(\bdx_*,\bdy_*)$ of an expanded mapping
\begin{equation}\label{dfm}
\mapform{\bdg}{\Sigma\subset\CC^m\times\CC^m}{\CC^n\times\CC^n\times\CC^{m-r}}{
(\bdx,\bdy)}{\left(\begin{array}{c}
\bdf(\bdx), \, J(\bdx)\,\bdy, \, R\,\bdy - \bde
\end{array}\right)}
\end{equation}
where $J(\bdx)$ is the Jacobian of $\bdf(\bdx)$, $R$ is a random 
$(m-r)\times m$ matrix and $\bde\,\ne\,\bdo$ with $r=\rank{J(\bdx_*)}$.
The deflation process terminates if $(\bdx_*,\bdy_*)$ is a regular zero
of $\bdg$ or, otherwise, continues recursively by expanding $\bdg$.
It is proved in \cite{DLZ,dayton-zeng} that the number of deflation steps is 
bounded by the {\em depth} of $\bdx_*$.
When depth-deflation terminates, the ultrasingular zero ~$\bdx_*$ 
of $\bdf$ is a component of the regular zero of the final expanded 
mapping.
As a result, the Gauss-Newton iteration locally converges to an accurate
zero at quadratic rate.
An earlier deflation strategy in \cite{lvz06} is also proven to terminate with
the number of steps bounded by the multiplicity.

By definition, a branch of $k$-dimensional semiregular zero of a mapping 
$\bdf$ can be parameterized as $\bdx=\phi(\bdz)$ for $\bdz$ in an open set.
As the parameter $\bdz$ varies, there is a significant likelihood that
$\nullity{\bdf_\bdx(\phi(\bdz))}$ degenerates below the dimension $k$ 
and reaches ultrasingularity. 
Such ultrasingular zeros can be of particular interest.
The following example shows that we can proactively seek such ultrasingularity
by applying the depth-deflation method.

\begin{example}[Ultrasingularity embedded in a semiregular solution set]
\em
The 
\newline cyclic-4 system arises in applications such as biunimodular
vectors that comes from a notion traces back to Gauss \cite{FuhRze}.
It is in the form of $\bdf(\bdx)\,=\,\bdo$ where 
$\bdx = (x_1,x_2,x_3,x_4)$ and
\begin{equation}\label{f4}  \bdf(\bdx) \,=\, \left[
\mbox{\scriptsize $\begin{array}{c}
x_1+x_2+x_3+x_4 \\
x_1\,x_2+x_2\,x_3+x_3\,x_4+x_4\,x_1 \\
x_1\,x_2\,x_3+x_2\,x_3\,x_4+x_3\,x_4\,x_1+x_4\,x_1\,x_2 \\
x_1\,x_2\,x_3\,x_4-1 \end{array}$}\right]
\end{equation}
The solution consists of two 1-dimensional branches
\begin{equation}\label{c4z}
 \{x_1\,=\,-x_3,\,x_2\,=\,-x_4,\,x_3\,x_4\,=\,\pm 1\}.
\end{equation}
All zeros in the branches are semiregular except eight ultrasingular zeros 
in the form
of $(\pm 1,\pm 1,\pm 1,\pm 1)$ ~and ~$(\pm i,\pm i,\pm i,\pm i)$ ~with 
proper choices of signs.
The cyclic-4 system becomes an ultrasingular equation at, say 
$\bdx_*=(1,-1,-1,1)$, for being 1-dimensional in (\ref{c4z}) but the nullity 
of the Jacobian $\bdf_\bdx(\bdx_*)$ is 2. 
However, it is a straightforward verification that, for almost all matrices 
$R\in\CC^{2\times 4}$, there is a unique $\bdy_*$ such that the point 
$(\bdx_*,\bdy_*)$ is a regular zero of the deflation mapping $\bdg$ 
in (\ref{dfm}).
As a result, the rank-$8$ Newton's iteration on $\bdg$ becomes the 
Gauss-Newton iteration that locally quadratically converges to 
$(\bdx_*,\bdy_*)$, solving the ultrasingular equation $\bdf(\bdx)=\bdo$.
The same assertion can be verified in the same way for all eight 
ultrasingular solutions.
The results show that, at least for cyclic-4 system the depth-deflation methods
deflates the ultrasingularity into regularity.
\end{example}

The rank-8 Newton's iteration on $\bdg$ converges specifically to those
eight ultrasingular zeros of $\bdf$ and  does {\em not} converges to 
other semiregular zeros in the same solution branch since they are not zeros 
of the deflation mapping $\bdg$ in (\ref{dfm}).
Consequently, the depth-deflation method can be proactively deployed to 
compute ultrasingular zeros if so desired.
At this point, however, the theories of the depth-deflation are lacking
at ultrasingularity embedded in semiregular branches of zeros and require
further studies.
Similar gaps exist in cases such as computing ultrasingular zeros at 
intersections of semiregular branches, and in cases where the entire 
branch of zeros are ultrasingular as shown in the following example
proposed by Barry Dayton.

\begin{example}[High dimension ultrasingularity]\label{e:hdu}
\em
~Consider the mapping below\linebreak with $\bdx = (x_1,\ldots,x_5)$
\begin{align*} 
\bdf(\bdx) \,=\, \left[\mbox{\scriptsize $\begin{array}{c}
x_2^2 x_4^2 + x_3^2 x_5^2 + x_1^3 - 2 x_2 x_4 \\
x_2^3 x_4^3 - 3 x_2^2 x_4^2 + x_3^2 x_5^2 + x_1^2 + 3 x_2 x_4 - 2 \\
x_3^3 x_5^3 + x_2^2 x_4^2 + x_1^2 - 2 x_2 x_4\end{array}$}\right]
\end{align*}
The solution set $S=\{(0,s,t,1/s,1/t) \,|\, s,t\ne 0\}$ is 2-dimensional but 
the nullity of the Jacobian is $4>2$, making the entire branch ultrasingular.
We apply the depth-deflation method by setting up the deflation mapping
$\bdg$ in (\ref{dfm}) with a random matrix $R\in\CC^{4\times 5}$.
For every $\bdx_*\in S$, there is a unique $\bdy_*$ such that 
$\bdg(\bdx_*,\bdy_*)=\bdo$.
Namely $\bdg$ also has a corresponding 2-dimensional zero set.
Anticipating this zero set to be semiregular, we set $r=10-2=8$ and
test the rank-8 Newton's iteration on $\bdg$ from an initial iterate
near $S$, say $\bdx_0=(\mbox{\scriptsize 0.001, .698, 1.201, 1.428, 0.833})$.
The rank-8 Newton's iteration converges to a point $(\tilde\bdx,\tilde\bdy)$
with the component $\tilde\bdx$ as
\[
\mbox{\scriptsize (.0, .699835056282962, 1.201681873936643, 
1.428908127739848, 0.832167000009791)}
\]
approximating a zero of $\bdf$ with an accuracy at hardware precision. 
The condition number $28.7$ indicates the Jacobian is indeed rank 8
and the solution is a semiregular zero of deflation mapping $\bdg$.

The result of this experiment shows that, at least for this polynomial system,
the depth-deflation method deflates the ultrasingularity into semiregularity.
\end{example}

Open questions remain such as: 
Does the depth-deflation deflates ultrasingularity in general? 
If so, under what conditions does the deflation terminate?
If not, are there proper modifications to overcome its limitations?
In fact, a numeric-symbolic deflation proposed by Hauenstein and Wampler
is proved to terminate in finitely many steps \cite{HauWam13}. 
Our preliminary experimental results also suggest the potential effectiveness
of the depth-deflation method combined with the novel low-rank Newton's 
iteration.

{\bf Acknowledgement.} We thank Dr. Wenrui Hao 
for the Bertini test 
and Dr. Tianran Chen for the HOM4PS test in Example~\ref{e:nag}.
We thank Dr. Barry Dayton for discussions on ultrasingular zeros 
and for suggesting Example~\ref{e:hdu}.
We also thank Dr. Jonathan Hauenstein for discussions on
deflation methods and results in \cite{HauWam13}.

\end{document}